\documentclass[12pt,reqno]{amsart}

\usepackage{amscd}
\usepackage{amssymb}
\usepackage[all]{xy}

\newtheorem{theorem}{Theorem}[section]
\newtheorem{proposition}[theorem]{Proposition}
\newtheorem{lemma}[theorem]{Lemma}

\theoremstyle{definition}
\newtheorem{definition}[theorem]{Definition}
\newtheorem{remark}[theorem]{Remark}
\newtheorem{Ass}[theorem]{Assumption}

\numberwithin{equation}{section}

\newcommand{\Aut}{\mathop{\rm Aut}\nolimits}
\newcommand{\res}{\mathop{\sf res}\nolimits}

\begin{document}

\title[Moduli spaces of framed logarithmic connections]{On the moduli spaces of framed logarithmic 
connections on a Riemann surface}

\author[I. Biswas]{Indranil Biswas}

\address{School of Mathematics, Tata Institute of Fundamental
Research, Homi Bhabha Road, Mumbai 400005, India}

\email{indranil@math.tifr.res.in}

\author[M. Inaba]{ Michi-aki Inaba}

\address{Department of Mathematics, Kyoto University, Kyoto 606-8502, Japan}

\email{inaba@math.kyoto-u.ac.jp}

\author[A. Komyo]{Arata Komyo}

\address{Center for Mathematical and Data Sciences, Kobe University,
1-1 Rokkodai-cho, Nada-ku, Kobe, 657-8501, Japan}

\email{akomyo@math.kobe-u.ac.jp}

\author[M.-H. Saito]{Masa-Hiko Saito}

\address{Department of Mathematics, Graduate School of Science,
Kobe University, Kobe, Rokko, 657-8501, Japan}

\email{mhsaito@math.kobe-u.ac.jp}



\date{}

\begin{abstract}
We describe some results on moduli space of logarithmic connections equipped
with framings on a $n$-pointed compact Riemann surface.
\end{abstract}

\maketitle

\section{Introduction}

Our aim here is to initiate the study of logarithmic connections, on a Riemann surface, equipped 
with a framing. To describe these objects, let $X$ be a compact connected Riemann surface and $D\, 
\subset\, X$ a finite subset. Given a holomorphic vector bundle $E$ on $X$ of rank $r$, a framing 
of it is an isomorphism $\phi_x\, :\, {\mathbb C}^r\, \longrightarrow\, E_x$ for every $x\, \in\, D$ 
(see \cite{BLP}). A framed logarithmic connection of rank $r$ on $(X,\, D)$ is a triple of the form 
$(E,\, \nabla,\, \phi)$, where $E$ is a holomorphic vector bundle on $X$ of rank $r$, $\nabla$ is a 
logarithmic connection on $E$ whose polar part is contained in $D$ and $\phi$ is a framing of $E$.

It can be shown that a moduli space of framed logarithmic connections exists as a Deligne-Mumford
stack. Furthermore, this moduli stack is equipped with a natural algebraic symplectic structure.
We note that the moduli space of logarithmic connections has a natural Poisson structure. The
forgetful map, from the moduli space of framed logarithmic connections to the moduli space
logarithmic connection that simply forgets the framing, is in fact Poisson.

Actually, the above results holds in a more general setting; see Section
\ref{section of symplectic form}. Also, the above results can be extended to
the set-up of holomorphic principal bundles with a logarithmic connection.

It is known that the moduli space of connections is not an affine variety. This should be 
contrasted with the fact that moduli space of connections is canonically biholomorphic to an 
affine scheme. Indeed, the character variety is affine, and the Riemann-Hilbert 
correspondence produces a biholomorphism of it with the moduli space of connections. This 
result is extended to the context of vector bundles equipped with parabolic
structure (see Theorem \ref{theorem on global algebraic function}).

The full details of the proofs will appear in \cite{BIKS}.

\section{Symplectic form on the moduli space of framed connections}\label{section of 
symplectic form}

Let $X$ be a compact connected Riemann surface of genus $g$,
and let $D$ be a reduced effective divisor on $X$. The holomorphic line bundle
$K_X\otimes {\mathcal O}_X(D)$, where $K_X$ is the holomorphic
cotangent bundle of $X$, will be denoted by $K_X(D)$.
Fix a connected complex algebraic proper subgroup
$H_x \,\subsetneq\, \mathrm{GL}(r,\mathbb{C})$
for each $x\in D$, and set $H\,=\, \{ H_x \}_{x \in D}$.

Let $E$ be a holomorphic vector bundle on $X$ of rank $r$. A {\it framing}
of $E$ is an ${\mathcal O}_D$--linear isomorphism
$\phi\,\colon\, E|_D \,\stackrel{\sim}{\longrightarrow}\,
{\mathcal O}_D^{\oplus r}$. Two framings $\phi$ and $\phi'$ are called equivalent if
$\phi'\circ\phi^{-1}|_x\,\in\, H_x$ for all $x\,\in\, D$.
A framing of $E$ with respect to $H$ is an equivalence class of framings of $E$.

Let $\mathfrak{h}_x$ be the Lie algebra of $H_x$.
The orthogonal complement of $\mathfrak{h}_x$,
with respect to the trace pairing $(A,\, B)\, \longmapsto\, \text{trace}(AB)$,
will be denoted by $\mathfrak{h}_x^{\bot}$.

Let $F(E)\, \longrightarrow\, X$ be the frame bundle for $E$; it is also the principal
$\text{GL}(r,{\mathbb C})$--bundle associated to $E$. Giving
a framing of $E$ with respect to $H$ is equivalent to giving a reduction of structure group
of $F(X)_x\, \longrightarrow\, x$ to $H_x$ for every $x\, \in\, D$. 
Given a framing of $E$ with respect to $H$, choose a framing $\phi$ in the equivalence class,
and consider the corresponding isomorphism $\text{Lie}(\text{GL}(r,{\mathbb C}))\,=\,
{\rm M}(r,{\mathbb C})\,\longrightarrow\, \text{End}(E)_x$ for every $x\,\in\, D$. The
restriction of this homomorphism to $\mathfrak{h}^{\bot}_x\, \subset\, {\rm M}(r,{\mathbb C})$ depends
on the choice of $\phi$ in the equivalence class, but the subspace
$\text{image}(\mathfrak{h}^{\bot}_x)\, \subset\,\text{End}(E)_x$ is independent of the choice
of $\phi$. This subspace of $\text{End}(E)_x$ will be denoted by $[\phi](\mathfrak{h}^{\bot}_x)$.
Similarly, the subspace
$\text{image}(\mathfrak{h}_x)\, \subset\,\text{End}(E)_x$ is independent of the choice
of $\phi$. This subspace of $\text{End}(E)_x$ will be denoted by $[\phi](\mathfrak{h}_x)$.

\begin{definition}
A framed connection with respect to $H$ is a triple
$(E,\,[\phi],\,\nabla)$, where
\begin{itemize}
\item $E$ is a holomorphic vector bundle on $X$ of rank $r$ equipped with
a framing $[\phi]$ with respect to $H$, and

\item $\nabla\,\colon\, E\,\longrightarrow \,E\otimes K_X(D)$ is a logarithmic connection
on $E$ such that $\res_x(\nabla)\,\in\, [\phi](\mathfrak{h}^{\bot}_x)$
for every $x\,\in\, D$.
\end{itemize}
\end{definition}

The group of all automorphisms $T$ of $E$ preserving $[\phi]$ such that
$(T\otimes {\rm Id}_{K_X(D)})\circ\nabla\,=\, \nabla\circ T$ will be denoted
by $\Aut(E,[\phi],\nabla)$.

\begin{definition}
A framed connection $(E,\,[\phi],\,\nabla)$ with respect to $H$ is {\it simple} if
the quotient group
$\Aut(E,[\phi],\nabla)\;/\; (\mathbb{C}^*\mathrm{id}_E\cap \Aut(E,[\phi],\nabla))$
is a finite group.
\end{definition}

Let $\mathcal{M}^{H}_{\rm{FC}}(d)$ denote
the moduli space of all framed connections
$(E,\,[\phi],\,\nabla)$ with respect to $H$ such that $\deg E\,=\,d$.

\begin{theorem}
The moduli space
$\mathcal{M}^{H}_{\rm{FC}}(d)$ exists as a Deligne-Mumford stack.
\end{theorem}

Given a framed connection $(E,\,[\phi],\,\nabla)$ with respect to $H$, we
construct a complex of sheaves ${\mathcal D}^{\bullet}$ as follows:
\begin{align*}
 {\mathcal D}^0\,
 &=\,
 \left\{ u\in {\mathcal End}(E) 
 \,\, \mid \,\, u(x) \,\in\, [\phi](\mathfrak{h}_x)
\,\ \forall \ x\,\in\, D \right\} 
 \\
 {\mathcal D}^1\,
 &=\,
 \left\{ v\in {\mathcal End}(E)\otimes K_X(D) 
 \,\,\mid\,\, \res_x(v)\,\in\, [\phi](\mathfrak{h}_x^{\bot}) \,\ \forall \ x\,\in\, D \right\} 
 \\
 {\mathcal D}^0 \ 
 &\ni \ u \ \longmapsto \ \nabla\circ u-(u\otimes\mathrm{id})\circ \nabla
 \ \in \ {\mathcal D}^1.
\end{align*}

\begin{theorem}
The tangent space of $\mathcal{M}^{H}_{\rm{FC}}(d)$ at any point
$(E,\,[\phi],\,\nabla)\,\in\, \mathcal{M}^{H}_{\rm{FC}}(d)$ is identified with 
the hypercohomology $\mathbb{H}^1({\mathcal D}^{\bullet})$ of the above
complex ${\mathcal D}^{\bullet}$.
\end{theorem}

Consider the following homomorphism of complexes
${\mathcal D}^{\bullet}\otimes{\mathcal D}^{\bullet}\, \longrightarrow\,
\Omega^{\bullet}_X$:
\begin{align*}
 ({\mathcal D}^{\bullet}\otimes{\mathcal D}^{\bullet})^0
 ={\mathcal D}^0\otimes{\mathcal D}^0
 \ \ni \ f\otimes g \ \longmapsto \ {\rm trace}(fg) \ \in {\mathcal O}_X
 \\
 ({\mathcal D}^{\bullet}\otimes{\mathcal D}^{\bullet})^1
 ={\mathcal D}^0\otimes{\mathcal D}^1\oplus
 {\mathcal D}^1\otimes{\mathcal D}^0
 \ \ni \ (f\otimes\omega,\eta\otimes g) \ \longmapsto \
 {\rm trace}(f\omega-g\eta) \ \in \ \Omega^1_X
\end{align*}
It produces a pairing
\begin{equation}\label{equation: pairing defining symplectic form}
 \Theta^H \,\colon\,
 \mathbb{H}^1({\mathcal D}^{\bullet})\otimes \mathbb{H}^1({\mathcal D}^{\bullet})
\,\longrightarrow\, \mathbb{H}^2({\mathcal D}^{\bullet}\otimes{\mathcal D}^{\bullet})
\, \longrightarrow\, \mathbf{H}^2(\Omega^{\bullet}_X)\,=\,\mathbb{C}\, .
\end{equation}

We have the following theorem.

\begin{theorem}\label{theorem of symplectic form}
The pairing $\Theta^H$ in \eqref{equation: pairing defining symplectic form} defines
a symplectic form on $\mathcal {M}^{H}_{\rm{FC}}(d)$.
\end{theorem}

We give an outline of the proof of Theorem \ref {theorem of symplectic form}.
We can check that
$\Theta^H$ is the skew-symmetric by a calculation in \u{C}ech cohomology.
The homomorphism
$\mathbb{H}^1({\mathcal D}^{\bullet})\,\longrightarrow\, \mathbb{H}^1({\mathcal D}^{\bullet})^*$
induced by $\Theta^H$ is isomorphic, because the vertical arrows except the middle one
in the exact commutative diagram
\[
\begin{CD}
H^0({\mathcal D}^0) @>>> H^0({\mathcal D}^1)
@>>> \mathbb{H}^1({\mathcal D}^{\bullet}) @>>> H^1({\mathcal D}^0)
@>>> H^1({\mathcal D}^1) \\
@VVV @VVV @VVV @VVV @VVV \\
H^1({\mathcal D}^1)^* @>>> H^1({\mathcal D}^0)^* @>>>
\mathbb{H}^1({\mathcal D}^{\bullet})^* @>>> H^0({\mathcal D}^1)^*
@>>> H^0({\mathcal D}^0)^* 
\end{CD}
\]
are Serre duality isomorphisms.
In other words, the $2$-form $\Theta^H$ is nondegenerate.

It remains to prove that the $2$-form $\Theta^H$ is $d$-closed.
Let us consider the moduli space
${\mathcal M}^e_{\rm FC}(d)$ of framed connections with respect to the trivial group.
When $H_x\,=\, \{e\}$, then ${\mathfrak h}^{\bot}_x\,=\, {\rm M}(r,{\mathbb C})$. Let
${\mathcal M}^{e}_{\rm FC}(d)^{\mathfrak{h}^{\bot}}\, \subset\,
{\mathcal M}^e_{\rm FC}(d)$ be the locus of all $(E,\, \phi,\, \nabla)$ such that
$\mathrm{res}_x(\nabla)$ lies in the image of $\mathfrak{h}_x^{\bot}$ for all $x\,\in\, D$.
Then we have the following diagram
\begin{equation*}
\xymatrix{
{\mathcal M}^{e}_{\rm FC}(d)^{\mathfrak{h}^{\bot}} \ \ar@ {^{(}->}[r]^-{\iota}
\ar[d]_-{\pi} &\mathcal{M}^e_{\rm{FC}} (d) \\
\mathcal{M}^H_{\rm{FC}} (d)\rlap{.}&
}
\end{equation*}
The above morphism $\pi$ sends any $(E,\,\phi,\,\nabla)$ to $(E,\,[\phi],\,\nabla)$,
and it is a smooth morphism.
So $\Theta^H$ is $d$-closed if $\pi^*\Theta^H$ is $d$-closed.
On the other hand, we can check the equality
$\pi^*\Theta^H=\iota^*\Theta^e$.
So it suffices to prove the theorem for the case $H\,=\,\{e\}$.
Since the moduli space ${\mathcal M}^e_{\rm FC}(d)$ is irreducible,
it is enough to prove that the restriction of 
$\Theta^e$ to some non-empty open subset of
${\mathcal M}^e_{\rm FC}(d)$ is $d$-closed.

For simplicity, we consider the case where $g\geq 2$.
Let us consider the following moduli spaces
\begin{align*}
 {\mathcal N}(d)
 &=\left\{ E \ \,\mid\, \ \text{$E$ is a stable bundle of rank $r$ and degree $d$} \right\}
 \\
 {\mathcal N}^e(d)
 &=
 \left\{ (E,\,\phi) \ \middle|\
 \begin{array}{l}
 \text{$E$ is a stable bundle of rank $r$ and degree $d$} \\
 \text{and $\phi\colon E|_D\xrightarrow{\sim} {\mathcal O}_D^{\oplus r}$ is an isomorphism}
 \end{array} \right\}
 \\
 {\mathcal M}(d)_0
 &=
 \left\{ (E,\,\nabla) \, \middle|\,
 \begin{array}{l}
 \text{$E$ is a stable bundle of rank $r$ and degree $d$,} \\
 \text{$\nabla\colon E\longrightarrow E\otimes K_X(D)$ is a logarithmic connection such that} \\
 \text{$\mathrm{res}_{x_i}(\nabla)\,=\,0$ for $1\leq i\leq n-1$
 and $\mathrm{res}_{x_n}(\nabla)\,=\,-\frac{d}{r}\mathrm{id}$}
 \end{array}
 \right\}
 \\
 {\mathcal M}^e_{\rm FC}(d)_0
 &=
 \left\{ (E,\,\phi,\,\nabla)\,\in\,{\mathcal M}^e_{FC}(d) \, \middle| \,
 \text{$(E,\,\phi)\,\in\, {\mathcal N}^e(d)$\, and\, $(E,\nabla)\,\in\, {\mathcal M}(d)_0$}
 \right\}.
\end{align*}
Then we have the following diagram
\begin{equation*}
\xymatrix{
{\mathcal M}(d)_0 \ar[d]^{p_0} &
\mathcal{M}^e_{\rm{FC}}(d )_0 \ar[l]_{\tilde{q}} \ar@ {^{(}->}[r]^-{\iota} \ar[d]_-{p^e_0} 
&\mathcal{M}^e_{\rm{FC}} (d) \ar[ld]^{p^e} \\
{\mathcal N}(d) \ar@/^8pt/@{.>}[u]^s & {\mathcal N}^e(d) \ar[l]_q \rlap{.}&
}
\end{equation*}
whose left square is a Cartesian diagram.
We can see that $p_0,\,p^e_0$ are smooth morphisms
and $\iota$ is a locally closed immersion.
So we can take a non-empty analytic open subset (or an \'etale neighborhood)
$U\,\subset\, {\mathcal N}(d)$ with a section
$s\,\colon\, U\,\longrightarrow\, p_0^{-1}(U)$ of $p_0$.
Let $U^e$ be the pullback of $U$ to ${\mathcal N}^e(d)$.
Then $s$ induces a section
\begin{equation}\label{equation: construction of section}
s^e\,\colon\, U^e\longrightarrow (p^e_0)^{-1}(U^e)
\,\hookrightarrow \,\mathcal{M}^e_{\rm{FC}} (d)
\end{equation}
of $p^e_0$.
Using the section $s^e$, we obtain the following isomorphism
$$
P_1\colon T^*U^e \,\xrightarrow{\ \cong \ }\, (p^e)^{-1}(U^e), \ \
(y,\,v) \,\longmapsto\, s^e(y) + v.
$$

\begin{lemma}\label{2019.11.19.15.22}
Let $\Phi_U$ be the Liouville $2$-form on the cotangent bundle $T^*U^e$.
Then,
\begin{equation*}
\Theta^e - (P_1^{-1})^* \Phi_U \,=\, (p^e)^* ((s^e)^* \Theta^e).
\end{equation*}
\end{lemma}

In view of Lemma \ref {2019.11.19.15.22}
it is enough to prove that $(s^e)^*\Theta^e$ is $d$-closed,
because the Liouville form is $d$-closed.
We can construct a two form
$\Theta_0$ on the moduli space
${\mathcal M}(d)_0$ in the same way as done in
\eqref{equation: pairing defining symplectic form}.
By construction, it is just the pullback of the Goldman symplectic form (\cite{Go},
\cite{AB}) via the Riemann-Hilbert morphism.
So we have $d\Theta_0\,=\,0$.

\begin{lemma}
The equality $\iota^*\Theta^e\,=\,\tilde{q}^*\Theta_0$ holds.
\end{lemma}

From the above lemma we can see that
$\iota^*\Theta^e$ is $d$-closed. Therefore,
so is $(s^e)^*\iota^*\Theta^e\,=\,(s^e)^*\Theta^e$,
and we are done.

\subsection{A Poisson map}

In this subsection we assume that $H_x\,=\,\{e\}$ for all $x\, \in\, D$.
Let $\mathcal{M}^{H}_{\rm{C}}(d)$ be the moduli space of pairs $(E,\, \nabla)$,
where $E$ is a holomorphic vector bundle on $X$ of rank $r$ and degree $d$, and $\nabla$
is a logarithmic connection on $E$ whose singular part is contained in $D$. Given any
$(E,\, \nabla)\, in\,\mathcal{M}^{H}_{\rm{C}}(d)$,
construct a complex $\widetilde{\mathcal D}^{\bullet}$ as follows:
$$
\widetilde{\mathcal D}^0\,=\, {\mathcal End}(E),\ \
\widetilde{\mathcal D}^1\,=\, {\mathcal End}(E)\otimes K_X(D)\, ,
$$
and $\widetilde{\mathcal D}^0\,\ni\, u \, \longmapsto \,
\nabla\circ u-(u\otimes\mathrm{id})\circ \nabla \, \in \, \widetilde{\mathcal D}^1$.
Then we have
$$
T_{(E,\nabla)}\mathcal{M}^{H}_{\rm{C}}(d) \,=\, {\mathbb H}^1(
\widetilde{\mathcal D}^{\bullet})\, .
$$
The homomorphism of hypercohomologies corresponding to
the natural inclusion of the Serre dual complex $\widetilde{\mathcal D}^{\bullet}$
in $\widetilde{\mathcal D}^{\bullet}$ produces a Poisson structure on the moduli
space $\mathcal{M}^{H}_{\rm{C}}(d)$.

For $(E,\, \nabla,\, \phi)\, \in\,\mathcal{M}^{H}_{\rm{FC}}(d)$, the differential
$T_{(E,\nabla,\phi)}\mathcal{M}^{H}_{\rm{FC}}(d)\, \longrightarrow\,
T_{(E,\nabla)}\mathcal{M}^{H}_{\rm{C}}(d)$ of the natural forgetful map
\begin{equation}\label{fm}
\mathcal{M}^{H}_{\rm{FC}}(d)\, \longrightarrow\,\mathcal{M}^{H}_{\rm{C}}(d)
\end{equation}
that forgets the framing coincides with the
homomorphism of hypercohomologies given by the inclusion map of the complex
${\mathcal D}^{\bullet}$ in $\widetilde{\mathcal D}^{\bullet}$. Using this it follows
that the forgetful map in \eqref{fm} is Poisson.

\section{Moduli space of parabolic connections}

It is known that there is no non-constant global algebraic function
on the moduli space of logarithmic connections
with central residues on a curve of genus at least $3$ \cite{BR}.
On the other hand, the character variety, which is a moduli space of
representations of a fundamental group, is affine.
So we can see that the Riemann-Hilbert morphism,
from the moduli space of connections
to the character variety, is not algebraic.
This non-algebraic map preserves the algebraic
symplectic forms on these two moduli spaces \cite{Bi}.

We replace the moduli space ${\mathcal M}^B_{\rm{FC}}$ with
the moduli space of the more general objects of parabolic connections.
Let $X$ be a compact Riemann surface
and $D$ be a reduced effective divisor on $X$.
We fix data
$\boldsymbol{\nu}\,=\, \{\nu^x_j\,\in\, \mathbb C\}_{1\leq j\leq r;\, x\in D}$ satisfying
\[
 \sum_{x\in D}\sum_{j=1}^r \nu^x_j\,=\, 0.
\]

\begin{definition}
We say that a triple
$(E,l,\nabla)$ is a $\boldsymbol{\nu}$-parabolic connection
if
\begin{itemize}
\item[(i)]
$E$ is a vector bundle of rank $r$ and degree $0$,
\item[(ii)]
$\nabla\colon E\longrightarrow E\otimes K_X(D)$
is a connection admitting poles along $D$
and
\item[(iii)]
$l$ is a filtration
$E|_D=l_1\supset l_2\supset\cdots\supset l_r\supset l_{r+1}=0$
satisfying the condition
$\left( \mathrm{res}_D(\nabla) - \nu_j\mathrm{id}\right) (l_j) \,\subset\, l_{j+1}$
for all $j\,=\,1,\,\cdots,\,r$.
\end{itemize}
\end{definition}

If we take $\boldsymbol{\nu}\,=\,\boldsymbol{0}$ so that all $\nu_j=0$,
then a $\boldsymbol{0}$-parabolic connection is equivalent to
a framed connection with respect to the Borel subgroup $B$.

For simplicity we adopt the following genericity assumption
on $\boldsymbol{\nu}$.

\begin{Ass}\label{Assmption of irreducible exponent}
For any integer $1\,\leq\, s\,<\,r$ and for any choice
of $s$ elements $\{j^x_1,\,\cdots,\,j^x_s\}$
in $\{1,\,\cdots,\,r\}$ for each $x\,\in \,D$, the following
holds:
\[
 \sum_{x\in D}\sum_{k=1}^s \nu_{j^x_k}\,\notin\,\mathbb{Z}.
\]
\end{Ass}

Let ${\mathcal M}(\boldsymbol{\nu})$
be the moduli space of $\boldsymbol{\nu}$-parabolic connections.
Assumption \ref{Assmption of irreducible exponent} ensures that
any $\boldsymbol{\nu}$-parabolic connection is irreducible.
So it is stable with respect to any parabolic weight.
Let ${\mathcal M}_{\rm{SL}}(\boldsymbol{\nu})$
be the moduli space of $\boldsymbol{\nu}$-parabolic connections
$(E,\,\nabla,\,l)$ with
$\det(E,\,\nabla)\,\cong\, ({\mathcal O}_X,\,d)$.
There is a closed immersion
$\iota\,\colon\, {\mathcal M}_{\rm{SL}}(\boldsymbol{\nu})
\,\hookrightarrow\, {\mathcal M}(\boldsymbol{\nu})$,
and also ${\mathcal M}_{\rm{SL}}(\boldsymbol{\nu})$ is smooth.

\begin{theorem}\label{theorem on global algebraic function}
Assume that $r\geq 2$, $n\geq 1$ and $g\geq 2$.
Then the moduli space ${\mathcal M}_{\rm{SL}}(\boldsymbol{\nu})$ 
of $\boldsymbol{\nu}$-connections
with the trivial determinant is not affine.
\end{theorem}

\begin{remark}\mbox{}
\begin{itemize}
\item[(1)]
There is a Riemann-Hilbert morphism from the moduli space of
$\boldsymbol{\nu}$-parabolic connections
to the character variety parameterizing the representations
of the fundamental group $\pi_1(X\setminus D)$
with a fixed local monodromy data. 
If we assume, in addition to Assumption \ref {Assmption of irreducible exponent},
that $\nu_j^x-\nu_k^x\,\notin\,\mathbb{Z}$ for $j\neq k$,
then the Riemann-Hilbert morphism is an analytic isomorphism.
Since the character variety is affine, 
the property that ${\mathcal M}_{\rm{SL}}(\boldsymbol{\nu})$ 
is not affine implies that the Riemann-Hilbert morphism is not algebraic.

\item[(2)] For a moduli space of framed connections with respect to the
Borel subgroup $B$ it is better to consider a stability with respect to a parabolic weight.
In that case, 
Assumption \ref {Assmption of irreducible exponent}
does not hold
nor the assumption in (1) above.
The moduli space of $\boldsymbol{\nu}$-parabolic connections
in the case of non-generic exponent $\boldsymbol{\nu}$, involving the case of $g=0$
is also interesting.
Indeed it contains loci of some special solutions of the isomonodromy equations,
which is like Riccatti loci in Painlev\'e VI equations
(see \cite {IIS1}, \cite{IIS2} and \cite {ST}).
The property of global algebraic functions on the moduli space of connections
in a special case of $g=0$
is also used in \cite{Arinkin}.

\item[(3)]
The formulation of the isomonodromy equation is given in
\cite{IIS1}, \cite{IIS2} and \cite{I1}, but it is also formulated in
\cite {BHH} and \cite{Ch}. Classically the isomonodromy equation is known
to be characterized by the isomonodromy $2$-form (\cite{JMU}, \cite {Iwasaki1}).
Since the isomonodromy equation gives an algebraic splitting of
the tangent bundle of the moduli space of connections,
we can construct the isomonodromy $2$-form from the symplectic form
constructed in Section \ref{section of symplectic form}.
\end{itemize}
\end{remark}

We will give an outline of the proof of Theorem \ref{theorem on global algebraic function},
using the following proposition.

\begin{proposition}\label{proposition: estimate of codimension} \
Under the assumption in Theorem \ref{theorem on global algebraic function},
consider the locus $Z$ in ${\mathcal M}_{\rm{SL}}(\boldsymbol{\nu})$
consisting of $\boldsymbol{\nu}$-parabolic connections
$(E,\,\nabla,\,l)$ whose underlying quasi-parabolic vector bundle $(E,\,l)$ is not simple.
Then the codimension of $Z$ in ${\mathcal M}_{\rm{SL}}(\boldsymbol{\nu})$
is at least $2$.
\end{proposition}

We denote by ${\mathcal N}^{\rm{spl}}_{\rm{par}}$
the moduli space of simple quasi-parabolic bundles
$(E,\,l)$ such that $\det E\,=\,{\mathcal O}_X$.
Let ${\mathcal M}^{\text{$\mathcal N$-spl}}_{\rm{SL}}(\boldsymbol{\nu})$
be the moduli space of $\boldsymbol{\nu}$-parabolic connections
with the trivial determinant such that
the underlying quasi-parabolic bundle $(E,\,l)$ is simple. Then there is a morphism
${\mathcal M}^{\text{$\mathcal N$-spl}}_{\rm{SL}}(\boldsymbol{\nu})
\longrightarrow {\mathcal N}^{\rm{spl}}_{\rm par}$
which is an affine space bundle. There is a universal family
$(\widetilde{E},\,\widetilde{l})$ on $X\times\mathcal{N}_{\rm par}^{\mathrm{spl}}$.
Let $\mathrm{At}(\widetilde{E})$ be the Atiyah bundle of $\widetilde{E}$ introduced in \cite{At},
which fits in the short exact sequence
\[
 0\, \longrightarrow \,{\mathcal End}(\widetilde{E})
 \,\longrightarrow\, \mathrm{At}(\widetilde{E})\,\longrightarrow\,
T_{X\times{\mathcal N}^{\mathrm{spl}}_{\rm par}/{\mathcal N}^{\mathrm{spl}}_{\rm par}}
 \,\longrightarrow\, 0.
\]
We set
\begin{align*}
 {\mathcal End}_{\mathfrak{sl}}(\widetilde{E},\,\widetilde{l})
 &:=
 \left\{ a\in {\mathcal End}(\widetilde{E})
 \, \middle| \,
 \text{$\mathrm{Tr}(a)=0$ and $a|_{D\times {\mathcal N}^{\mathrm{spl}}_{\rm par}}(\widetilde{l}_j)
 \subset \widetilde{l}_j$
 for any $j$} \right\}\\
 {\mathcal End}_{\mathfrak{sl}}^{\mathrm{nil}}(\widetilde{E},\,\widetilde{l})
 &=
 \left\{ a\in {\mathcal End}(\widetilde{E})
 \, \middle| \,
 \text{$\mathrm{Tr}(a)=0$ and $a|_{D\times {\mathcal N}^{\mathrm{spl}}_{\rm par}}(\widetilde{l}_j)
 \,\subset\, \widetilde{l}_{j+1}$
 for any $j$} \right\}.
\end{align*}
We can define a subsheaf
$\mathrm{At}(\widetilde{E},\widetilde{l})
\,\subset\, \mathrm{At}(\widetilde{E})$
fitting in the exact commutative diagram
\[
 \begin{CD}
 0 \, \longrightarrow \ & {\mathcal End}_{\mathfrak{sl}}(\widetilde{E},\,\widetilde{l})
 @>>>\mathrm{At}(\widetilde{E},\,\widetilde{l})
 @>>> T_{X\times{\mathcal N}^{\mathrm{spl}}_{\rm par}/{\mathcal N}^{\mathrm{spl}}_{\rm par}}
 (-D\times {\mathcal N}^{\mathrm{spl}}_{\rm par})
 & \ \longrightarrow \, 0
 \\
 & @VVV @VVV @VVV & \\
 0\, \longrightarrow \ & {\mathcal End}(\widetilde{E})
 @>>> \mathrm{At}(\widetilde{E}) @>>>
 T_{X\times{\mathcal N}^{\mathrm{spl}}_{\rm par}/{\mathcal N}^{\mathrm{spl}}_{\rm par}} 
 & \ \longrightarrow\, 0.
 \end{CD}
\]
Tensoring $K_X(D)$, we get an exact sequence
\[
 0 \longrightarrow 
 {\mathcal End}_{\mathfrak{sl}}(\widetilde{E},\,\widetilde{l}) \otimes K_X(D)
 \longrightarrow
 \mathrm{At}(\widetilde{E},\,\widetilde{l})\otimes K_X(D)
 \longrightarrow{\mathcal O}_{X\times{\mathcal N}^{\mathrm{spl}}_{\rm par}}
 \longrightarrow 0,
\]
from which
we get a short exact sequence
\[
 0 \longrightarrow 
 \pi_* \left( {\mathcal End}_{\mathfrak{sl}}(\widetilde{E},\,\widetilde{l})
 \otimes K_X(D) \right)
 \longrightarrow
 \pi_*\left(\mathrm{At}(\widetilde{E},\,\widetilde{l})\otimes K_X(D)\right)
 \longrightarrow\pi_*\left( {\mathcal O}_{X\times{\mathcal N}^{\mathrm{spl}}_{\rm par}}\right)
 \longrightarrow 0.
\]
We put
${\mathcal Q}\,:=\,\pi_*\left(\mathrm{At}(\widetilde{E},\,\widetilde{l})\otimes K_X(D)\right)$,
and consider the projective bundle
\[
 \mathbb{P}_*({\mathcal Q})\,=\,
 \mathrm{Proj}
 \left(
 \mathrm{Sym}^*
 \left( {\mathcal Q}^{\vee} \right)
 \right)
\]
where
$\mathrm{Sym}^*({\mathcal Q}^{\vee})$ is the symmetric algebra of
${\mathcal Q}^{\vee}$
over
${\mathcal O}_{{\mathcal N}^{\mathrm{spl}}_{\rm par}}$.
There is a tautological sub line bundle
\[
 {\mathcal O}_{\mathbb{P}_*({\mathcal Q})}(-1)\, \hookrightarrow\,
 Q\otimes{\mathcal O}_{\mathbb{P}_*({\mathcal Q})}.
\]
There are induced sections
\begin{align*}
 \epsilon
 &:
 {\mathcal O}_{\mathbb{P}_*({\mathcal Q})}(-1)
 \hookrightarrow
 Q\otimes{\mathcal O}_{\mathbb{P}_*({\mathcal Q})}
 \longrightarrow
 \pi_*\left( {\mathcal O}_{X\times{\mathcal N}^{\mathrm{spl}}_{\rm par}}\right)
 \otimes{\mathcal O}_{\mathbb{P}_*({\mathcal Q})}
 ={\mathcal O}_{\mathbb{P}_*({\mathcal Q})}
 \\
 \widetilde{\nu}_j
 &\colon
 {\mathcal O}_{\mathbb{P}_*({\mathcal Q})}(-1)
 \hookrightarrow
 Q\otimes{\mathcal O}_{\mathbb{P}_*({\mathcal Q})}
 \xrightarrow{\mathrm{res}_D}
 {\mathcal End}_{\mathfrak{sl}}(\widetilde{E},\widetilde{l})\big|_{D\times \mathbb{P}_*({\mathcal Q})}
 \longrightarrow
 {\mathcal End}(\widetilde{l}_j/\widetilde{l}_{j+1})\otimes{\mathcal O}_{\mathbb{P}_*({\mathcal Q})}
 ={\mathcal O}_{D\times\mathbb{P}_*({\mathcal Q})}
\end{align*}
Let
${\mathcal I}$ be the ideal sheaf of the graded
${\mathcal O}_{{\mathcal N}^{\mathrm{spl}}_{\rm par}}$-algebra
$\mathrm{Sym}^*({\mathcal Q}^{\vee})$
generated by
$$\left\{ \widetilde{\nu}_j|_x-\nu_j|_x\,\epsilon \: \middle| \ x\in D, 1\leq j\leq r \right\}\, .$$
Then there is a short exact sequence
\begin{equation}\label{equation: exact sequence for graded algebra}
0\, \longrightarrow\, \mathrm{Sym}^*( {\mathcal Q}^{\vee} )/{\mathcal I} 
\, \stackrel{\epsilon}\longrightarrow\, \mathrm{Sym}^*( {\mathcal Q}^{\vee} )/{\mathcal I}
\, \longrightarrow\, \mathrm{Sym}^*\left( 
\pi_* \left( {\mathcal End}_{\mathfrak{sl}}^{\mathrm{nil}}(\widetilde{E},\widetilde{l})
\otimes K_X(D) \right) \right)
\, \longrightarrow\, 0.
\end{equation}
We put
\[
 \overline {{\mathcal M'}}
 := \mathrm{Proj} (\mathrm{Sym}^*( {\mathcal Q}^{\vee} )/{\mathcal I} )
 \subset
 \mathbb{P}_*({\mathcal Q}).
\]
Let $Y\,\subset\, \overline {\mathcal M'}$ be
the effective divisor defined by $\epsilon\,=\,0$.
There is a canonical isomorphism
${\mathcal M}^{\text{$\mathcal N$-spl}}_{\rm{SL}}(\boldsymbol{\nu})
\,\cong\,\overline { \mathcal M' } \setminus Y\,\cong\,
\mathrm{Spec} \left( \mathrm{Sym}^*( {\mathcal Q}^{\vee} )/{\mathcal I} \right)_{(\epsilon)}
$, where $\left( \mathrm{Sym}^*( {\mathcal Q}^{\vee} )/{\mathcal I} \right)_{(\epsilon)}$
is the degree zero component of the localized graded ring
$\left( \mathrm{Sym}^*( {\mathcal Q}^{\vee} )/{\mathcal I} \right)_{\epsilon}$.

Suppose that 
${\mathcal M}_{\rm{SL}}(\boldsymbol{\nu})$
is affine.
Then the ring of global algebraic functions
\[
 A\,=\,
 \Gamma \left(
 {\mathcal M}_{\rm{SL}}(\boldsymbol{\nu}) ,\
 {\mathcal O}_{ {\mathcal M}_{\rm{SL}}(\boldsymbol{\nu}) }\right)
\]
is a finitely generated $\mathbb{C}$-algebra of Krull dimension
$2(r^2-1)(g-1)+r(r-1)n$.
By Proposition \ref{proposition: estimate of codimension},
\[
 A\,=\, 
 \Gamma\left({\mathcal M}^{\text{$\mathcal N$-spl}}_{\rm{SL}}(\boldsymbol{\nu}),\,
 {\mathcal O}_{ {\mathcal M}^{\text{$\mathcal N$-spl}}_{\rm{SL}}(\boldsymbol{\nu}) }\right)
\,=\,
 \Gamma\left ( {\mathcal N}^{\mathrm{spl}}_{\rm par},\,
 \left( \mathrm{Sym}^*( {\mathcal Q}^{\vee} )/{\mathcal I} \right)_{(\epsilon)} \right)
\,=\,
 \Gamma\left ( {\mathcal N}^{\mathrm{spl}}_{\rm par},\,
 \mathrm{Sym}^*( {\mathcal Q}^{\vee} )/{\mathcal I} \right)_{(\epsilon)}.
\]
Since $A$ is a finitely generated $\mathbb{C}$-algebra,
it can be proved that
$\Gamma\left ( {\mathcal N}^{\mathrm{spl}}_{\rm par},\,
 \mathrm{Sym}^*( {\mathcal Q}^{\vee} )/{\mathcal I} \right)$
is a finitely generated graded $\mathbb{C}$-algebra.
So the Krull dimension of
$\Gamma\left ( {\mathcal N}^{\mathrm{spl}}_{\rm par}, \,
\mathrm{Sym}^*( {\mathcal Q}^{\vee} )/{\mathcal I} \right)$
should be
$\text{Krull-dim}(A)+1$.
From the exact sequence (\ref{equation: exact sequence for graded algebra}),
we get the exact sequence
\begin{align*}
0 \,\longrightarrow\, 
\Gamma\left({\mathcal N}^{\rm{spl}}_{\rm par},\,
\mathrm{Sym}^*( {\mathcal Q}^{\vee} )/{\mathcal I} \right)
 &\,\stackrel{\epsilon}\longrightarrow\,
 \Gamma\left ( {\mathcal N}^{\mathrm{spl}}_{\rm par},\,
 \mathrm{Sym}^*( {\mathcal Q}^{\vee} )/{\mathcal I} \right)
 \\
 &\longrightarrow\,
 \Gamma \left( {\mathcal N}^{\mathrm{spl}}_{\rm par},\, 
 \mathrm{Sym}^*\left( 
 \pi_* \left( {\mathcal End}_{\mathfrak{sl}}^{\mathrm{nil}}(\widetilde{E},\,\widetilde{l})
\otimes \Omega^1_X(D) \right)\right)\right).
\end{align*}
So we have
\begin{align}\label{equation: inequality for Krull dimension}
\begin{split}
 \text{Krull-dim}(A)
 &\,=\,
 \text{Krull-dim}\left(
 \Gamma\left ( {\mathcal N}^{\mathrm{spl}}_{\rm par}, \,
 \mathrm{Sym}^*( {\mathcal Q}^{\vee} )/{\mathcal I} \right) \big/
 \epsilon \; \Gamma\left ( {\mathcal N}^{\mathrm{spl}}_{\rm par},\, 
 \mathrm{Sym}^*( {\mathcal Q}^{\vee} )/{\mathcal I} \right)
 \right)
 \\
 &=\,
 \mathrm{tr.deg}_{\mathbb{C}}
 \left(
 \Gamma\left ( {\mathcal N}^{\mathrm{spl}}_{\rm par}, \,
 \mathrm{Sym}^*( {\mathcal Q}^{\vee} )/{\mathcal I} \right) \big/
 \epsilon \; \Gamma\left ( {\mathcal N}^{\mathrm{spl}}_{\rm par},\, 
 \mathrm{Sym}^*( {\mathcal Q}^{\vee} )/{\mathcal I} \right)
 \right)
 \\
 &\leq\,
 \mathrm{tr.deg}_{\mathbb{C}}
 \Gamma \left( {\mathcal N}^{\mathrm{spl}}_{\rm par}, \,
 \mathrm{Sym}^*\left( 
 \pi_* \left( {\mathcal End}_{\mathfrak{sl}}^{\mathrm{nil}}(\widetilde{E},\,\widetilde{l})
 \otimes \Omega^1_X(D) \right)\right)\right) 
\end{split}
\end{align}

On the other hand, let ${\mathcal M}^{\rm{spl}}_{\mathrm{SL},\mathrm{Higgs}}$
be the moduli space of simple parabolic Higgs bundles
$(E,\,\Phi,\,l)$ such that $\Phi|_D(l_j)\subset l_{j+1}$ for $1\leq j\leq r$,
$\det(E)\,=\,{\mathcal O}_X$ and ${\rm trace}(\Phi)\,=\,0$.
It contains the moduli space
${\mathcal M}^{\boldsymbol{\alpha}}_{\mathrm{SL},\mathrm{Higgs}}$
of $\boldsymbol{\alpha}$-stable parabolic Higgs bundles
for any choice of parabolic weight $\boldsymbol{\alpha}$. Let
${\mathcal M}^{\text{$\mathcal N$-spl}}_{\mathrm{SL},\mathrm{Higgs}}$
be the open subspace of
${\mathcal M}^{\rm{spl}}_{\mathrm{SL},\mathrm{Higgs}}$
consisting of those whose underlying quasiparabolic bundle
$(E,\,l)$ is simple. Then ${\mathcal M}^{\text{$\mathcal N$-spl}}_{\mathrm{SL},\mathrm{Higgs}}$
is isomorphic to the cotangent bundle of
${\mathcal N}^{\mathrm{spl}}_{\rm par}$, which means that
\[
{\mathcal M}^{\text{$\mathcal N$-spl}}_{\mathrm{SL},\mathrm{Higgs}}\, \cong\,
 \mathrm{Spec} \left(
 \mathrm{Sym}^* \left(
 \pi_* \left( {\mathcal End}_{\mathfrak{sl}}^{\mathrm{nil}}(\widetilde{E},\,\widetilde{l})
 \otimes \Omega^1_X(D) \right)
\right) \right).
\]
We have the codimension estimation for the moduli space of parabolic Higgs bundles
which is similar to Proposition \ref {proposition: estimate of codimension}.
So we have
\begin{align}\label{equation: algebraic function of moduli of Higgs bundles}
\begin{split}
 \Gamma \left( {\mathcal N}^{\mathrm{spl}}_{\rm par}, \,
 \mathrm{Sym}^*\left( 
 \pi_* \left( {\mathcal End}_{\mathfrak{sl}}^{\mathrm{nil}}(\widetilde{E},\,\widetilde{l})
 \otimes \Omega^1_X(D) \right)\right)\right) 
 &\,\cong\,
 \Gamma \left(
 {\mathcal M}^{\text{$\mathcal N$-spl}}_{\mathrm{SL},\mathrm{Higgs}} ,\,
 {\mathcal O}_{{\mathcal M}^{\text{$\mathcal N$-spl}}_{\mathrm{SL},\mathrm{Higgs}} } \right)
 \\
 &=\,
 \Gamma \left( {\mathcal M}^{\rm{spl}}_{\mathrm{SL},\mathrm{Higgs}} ,\,
 {\mathcal O}_{ {\mathcal M}^{\rm{spl}}_{\mathrm{SL},\mathrm{Higgs}} } \right)
 \\
 &\subseteq\,
 \Gamma \left(
 {\mathcal M}^{\boldsymbol{\alpha}}_{\mathrm{SL},\mathrm{Higgs}} ,\,
 {\mathcal O}_{ {\mathcal M}^{\boldsymbol{\alpha}}_{\mathrm{SL},\mathrm{Higgs}} } \right).
\end{split}
\end{align}
Choosing $\boldsymbol{\alpha}$ generically, we may assume that
the $\boldsymbol{\alpha}$-semistability implies $\boldsymbol{\alpha}$-stability.
So the Hitchin map
\[
 {\mathcal M}^{\boldsymbol{\alpha}}_{\mathrm{SL},\mathrm{Higgs}}
\, \longrightarrow\, \bigoplus_{k=2}^r
 H^0\left(X, \, K_X^{\otimes k}((k-1)D)\right)
\]
is a proper morphism.
Hence the $\mathbb{C}$-algebra
$\Gamma \left(
 {\mathcal M}^{\boldsymbol{\alpha}}_{\mathrm{SL},\mathrm{Higgs}} ,\,
 {\mathcal O}_{ {\mathcal M}^{\boldsymbol{\alpha}}_{\mathrm{SL},\mathrm{Higgs}} } \right)$
is a finite algebra over
\[
 \Gamma \left(
 \bigoplus_{k=2}^r
 H^0\left(X, \, K_X^{\otimes k}((k-1)D)\right),
 {\mathcal O}_{
 \bigoplus_{k=2}^r
 H^0\left(X, K_X^{\otimes k}((k-1)D)\right)} \right),
\]
whose transcendence degree is the same as its Krull dimension
$r^2(g-1)-g+1+nr(r-1)/2$.
Combining \eqref{equation: inequality for Krull dimension}
and \eqref{equation: algebraic function of moduli of Higgs bundles}
we get the inequality
\begin{align*}
 2(r^2-1)(g-1)+r(r-1)n\,=\,\text{Krull-dim}(A)
 &\,\leq\,
 \mathrm{tr.deg}_{\mathbb{C}}
 \Gamma \left(
 {\mathcal M}^{\boldsymbol{\alpha}}_{\mathrm{SL},\mathrm{Higgs}} ,\,
 {\mathcal O}_{ {\mathcal M}^{\boldsymbol{\alpha}}_{\mathrm{SL},\mathrm{Higgs}} } \right)
 \\
 &\leq\,
(r^2-1)(g-1)+\frac{nr(r-1)} {2},
\end{align*}
which is a contradiction.


\end{document}